\documentclass{article}
\author{Claire Levaillant}

\title{Congruences related to Miki's identity}

\usepackage{amsmath, amssymb, amsthm, enumerate}

\newcommand{\pdn}{\left[\begin{array}{l}p\\2n\end{array}\right]}
\newcommand{\mpt}{\;\text{mod}\,p^3}
\newcommand{\mpd}{\;\text{mod}\,p^2}
\newcommand{\mpu}{\;\text{mod}\,p}
\newcommand{\dmz}{\delta_0(m_0)}
\newcommand{\di}{\delta_0(i)}
\newcommand{\dui}{\delta_1(i)}
\newcommand{\dk}{\delta_0(k)}

\newcommand{\dmzu}{\delta_1(m_0)}
\newcommand{\ms}{\mathcal{S}^{'}_}
\newcommand{\mc}{\mathcal{C}}
\newcommand{\mcu}{\mc_1}
\newcommand{\mcd}{\mc_2}
\newcommand{\mct}{\mc_3}
\newcommand{\mcq}{\mc_4}
\newcommand{\mcf}{\mc_5}
\newcommand{\mcs}{\mc_6}
\newcommand{\mb}{\mathcal{B}}
\newcommand{\si}{\left[\begin{array}{l}\;p\;\\2n+1\end{array}\right]}
\newcommand{\md}{\mathcal{D}}
\newcommand{\mh}{\mathcal{H}}
\newcommand{\bn}{\binom}
\newcommand{\nts}{\negthickspace}
\newcommand{\mt}{\mathcal{T}}
\newcommand{\as}{A^{\star}}
\newcommand{\ma}{\mathcal{A}}
\newcommand{\mcp}{\mathcal{P}}
\newcommand{\mpq}{\;\;\text{mod}\,p^4}

\begin{document}
\maketitle
\begin{center}
\textit{This paper is dedicated to Jennifer Whitehead}
$$\begin{array}{l}\\\end{array}$$
\end{center}
\begin{center}\textbf{Abstract.}\end{center}

Given an odd prime number $p$, we present three different ways, of comparable levels of elegance and sophistication, to relate modulo $p$ certain truncated convolutions of divided Bernoulli numbers to certain full convolutions of divided Bernoulli numbers.

One approach consists of computing modulo $p^3$ the Stirling numbers on $p$ letters with odd indices in two different ways: on one hand using some $p$-adic analysis on the polynomial $X^{p-1}+(p-1)!$ with $p$-adic integer coefficients (as in \cite{LEV} $2019$) -- where the truncated convolutions of divided Bernoulli numbers emerge -- and on the other hand using Newton's formulas -- where the full convolutions of divided Bernoulli numbers arise.

Another approach is based on Hiroo Miki's identity (\cite{MI} $1978$) relating a convolution of divided Bernoulli numbers to a binomial convolution of divided Bernoulli numbers and harmonic numbers. It also crucially uses an identity which computes an alternating binomial sum of harmonic numbers and which is due to Michael Spivey (\cite{SP} $2007$).

A third approach is based on computing modulo $p^3$ the multiple harmonic sums on the first $p-1$ integers, using both Newton's formulas and Zhi-Hong Sun's pioneering result (\cite{SU2} $2000$) for computing modulo $p^3$ the generalized harmonic numbers with even indices. These multiple harmonic sums get related to the Stirling numbers modulo $p^3$ by using James Whitbread Lee Glaisher's formula (\cite{GL} $1900$), namely $(p-1)!=p\,B_{p-1}-p\;\mpd$.

If the Stirling numbers with even indices are known modulo $p^3$ since Glaisher, little is known about the Stirling numbers with odd indices modulo $p^3$. Equivalently, we do not know how to resolve the convolution of divided Bernoulli numbers modulo $p$ in the general case. The situation is similar with multiple harmonic sums. Equivalently, we do not know how to resolve the truncated convolutions of divided Bernoulli numbers. Of course, we have proven that both problems are equivalent. However, a few special cases can be worked out.

For instance, Jianqiang Zhao (\cite{ZO} $2007$) determines the multiple harmonic sum on two integers modulo $p^3$. He relates it to the harmonic number $\mathcal{H}_{p-1}$ on $p-1$ integers. We extend by our method Zhao's result to the multiple harmonic sum on four integers. Moreover, the full convolution of divided Bernoulli numbers of order $p-5$ corresponds to the case when there is a single term in the truncated convolution and is thus revealed. From there, we also deduce the Stirling number on five disjoint cycles modulo $p^3$. 

The implication of Zhao's result with the surrounding works goes as follows. 
Wolstenholme's theorem asserts that $\mh_{p-1}$ is divisible by $p^2$. Sun provides a way for computing the generalized harmonic numbers to the modulus $p^3$. By combining our method, Sun's work applied to $H_2$ and Zhao's result, we are able to provide $\mh_{p-1}$ modulo $p^4$. This is a result which pre-existed in Sun's work without proof.

Additionally, from knowing the multiple harmonic sum on two integers modulo $p^3$ and $\mathcal{H}_{p-1}$ modulo $p^4$, we derive the Stirling number on three cycles modulo $p^3$ and then deduce the convolution of divided Bernoulli numbers of order $p-3$ modulo $p$.

As part of our work, we also obtain the inverse of $(p-1)!$ modulo $p^3$, namely we have the following congruence in terms of divided Bernoulli numbers.
$$\frac{1}{(p-1)!}=3(p\mb_{p-1}-1)-p\,\mb_{2(p-1)}-\frac{1}{2}p^2\mb_{p-1}^2\;\;\mpt$$
This inverse is the multiple harmonic sum on $p-1$ integers. Our proof relies on a generalization of Wilson's theorem to the modulus $p^3$ by Zhi-Hong Sun as well as his computation for the generalized harmonic number $H_{p-1}$ modulo $p^3$.

Interestingly, our work shows that if $p$ is a Wilson prime and $\mb_{p-1+(p-1-2n)}=\mb_{p-1-2n}\;\mpd$ for some $2n\in\lbrace 4,6,\dots,p-7\rbrace$, then the truncated convolution $\mt\mc\mb(p+1-2n,p-3)$ is the opposite of the full convolution of order $p-1-2n$ modulo $p$.\\ It is unknown whether or not there would exist such primes. \\
It was shown by Richard Crandall, Karl Dilcher and Carl Pomerance that there are only three Wilson primes below $5\times 10^8$, namely $5$, $13$ and $563$ (\cite{CDP} $1997$).\\
As for the first congruence above on the divided Bernoulli numbers, it holds modulo $p$ as it is a congruence of K\"ummer type. \\
Reijo Ernvall and Tauno Mets\"ankyl\"a have shown that if $(p,p-1-2n)$ is an irregular pair, then the congruence cannot hold modulo $p^2$. \\
When $p=13$, the choices for $2n$ are $4$ or $6$ and the congruences are not satisfied modulo $p^2$:
$\mb_{20}\not\equiv\mb_8\;\text{mod}\,13^2$ and $\mb_{18}\not\equiv\mb_6\;\text{mod}\,13^2$.

\newpage
\begin{center}\huge{\textbf{Contents}}
\end{center}
$$\begin{array}{l}\\\\\end{array}$$
\textbf{1. Introduction, notations and useful results}\dotfill$\mathbf{4}$\\

1.1. Introduction and main results\dotfill $4$ \\

1.2. Notations and useful results\dotfill$9$\\\\
\textbf{2. The Stirling numbers modulo $p^3$}\dotfill$\mathbf{13}$\\

2.1. Stirling numbers with even indices\dotfill$13$\\

2.2. Stirling numebrs with odd indices\dotfill$17$\\

\hspace{0.25in}2.2.1. General calculation\dotfill$18$\\

\hspace{0.25in}2.2.2. Combinatorial analysis of the coefficients $D_s$'s and $G_k$'s\dotfill$21$ \\

\hspace{0.25in}2.2.3. Where we resolve the truncated convolutions\dotfill$22$ \\

\hspace{0.25in}2.2.4. Final expression and confrontation with the result issued from Manner I\dotfill$30$\\

\noindent \textbf{3. The multiple harmonic sums modulo $p^3$}\dotfill$\mathbf{32}$ \\\\
\textbf{4. Concluding words}\dotfill$\mathbf{34}$

\newpage
\section{Introduction, notations and useful results}
\subsection{Introduction and Main results}
\indent Unless otherwise mentioned, $p$ will be a given odd prime. By times it will be necessary to impose $p\geq 5$ or even $p\geq 7$. But our work is interesting for large primes as we deal with symmetric groups of order $p!$ and small orders can be easily worked out by hand. \\
This paper arose from the interest in computing the sum $(28)$ of \cite{LEV} modulo $p$, namely
$$\sum_{r=1}^{k-1}\frac{B_{2r}B_{2k-2r}}{2r}\;\text{mod}\,p,$$
in order to obtain a nice expression for the Stirling numbers $\left[\begin{array}{l}\;\;\;p\\p-2k\end{array}\right]$ modulo $p^3$.
Stirling numbers, as just denoted by their brackets, with odd or even indices $s$ are interesting numbers because, as one of their combinatorial facets, they count the number of permutations on $p$ letters that decompose into a disjoint product of $s$ cycles. \\Knowing the $p$-adic expansion of these numbers can be useful in a variety of situations. Earlier on, British mathematician Glaisher was already aware of that. Around the turn of the $20$th century, he was first to realize their $p$-adic analysis up to the modulus $p^2$ and up to the modulus $p^3$ for the even indices. The reference on the topic is \cite{GL}, a work that is also based on \cite{GL1}.\\
A method based on Newton's formulas and originating a century later in Sun's work \cite{SU2} had allowed us in \cite{LEV} to:
$$\begin{array}{l}
\text{(i) Retrieve Glaisher's formula for the even indices, see \cite{LEV}.}\\
\text{(ii) Obtain the congruence $(28)$ of \cite{LEV} for the odd indices.}
\end{array}$$
For the details regarding points $(i)$ and $(ii)$, the reader is referred to forthcoming Result $8$ of the current introduction which states the result of \cite{LEV}. \\
For convenience, we will refer to this method as "Manner I".\\

\noindent Our original goal presently was to compute the Stirling numbers modulo $p^3$ with odd indices in a different way, say "Manner II", namely one which uses some elementary $p$-adic analysis on the polynomial $X^{p-1}+(p-1)!$ with $p$-adic integer coefficients. Unfortunately, after far more efforts than those deployed for Manner I, we obtain the same formula as previously, thus failing at resolving the convolution. However, in order to show the certainly non-trivial equivalence between the results provided independently by Manner I and Manner II, we expand on resolving congruences involving truncated convolutions of divided Bernoulli numbers in terms of full convolutions of such numbers, using Hiroo Miki's identity. These are interesting congruences in their own right.
Miki's identity relates convolutions of divided Bernoulli numbers to binomial convolutions of divided Bernoulli numbers and harmonic numbers.\\
Another approach to our work is thus the following. The congruence that we get between the two types of convolutions which is derived from a central use of Miki's identity and which we just discussed above can be obtained from confronting Manner I and Manner II for the computation of the Stirling numbers on $p$ letters with odd indices modulo $p^3$, independently from Miki's identity. We note the strong interconnection between Miki's identity and the Stirling numbers.
Ira Gessel had already highlighted this interconnection in \cite{GE}, when he gave a simpler proof than Miki's original proof for Miki's identity. However, he uses Stirling numbers of the second kind instead of the first kind, that is the number of ways to partition a set of $n$ objects into $k$ non-empty subsets (denoted by $\left\lbrace\begin{array}{l}n\\k\end{array}\right\rbrace$). Gessel's proof is based on two different expressions for Stirling numbers of the second kind. We note that Miki's original proof was based on some p-adic analysis and also used the Fermat quotient $q_p(a)$. Another proof of Miki's identity using $p$-adic analysis appears in \cite{SHY}. So far we have seen two approaches, namely the "Miki approach" and the "polynomial approach" for deriving congruences concerning these truncated convolutions of divided Bernoulli numbers.\\
We now describe a third versant of this research. By using Newton's formulas, we are able to relate the multiple harmonic sums modulo $p^3$ to the truncated convolutions. Then, we are able to relate the multiple harmonic sums $\as_{2n}$ (sums of products of $2n$ reciprocals amongst $p-1$) to the Stirling numbers $A_{p-1-2n}$ onto $(2n+1)$ cycles. Moreover, when $2n=p-1$, the truncated convolution is the full convolution which is also related to $(p-1)!$. Then, we obtain $\as_{p-1}$ modulo $p^3$, that is we obtain the inverse of $(p-1)!$ modulo $p^3$.\\
The method also provides a third and independent way of relating the truncated convolutions to the full convolutions modulo $p^3$. \\Since the first Stirling numbers $A_4$, $A_6$, $A_8$, etc may be computed easily modulo $p^3$ using the values of the Bernoulli numbers with small indices, we obtain some neat formulas for the mirror multiple harmonic sums $\as_{p-5}$, $\as_{p-7}$, $\as_{p-9}$, etc.


Establishing congruences modulo $p^r$ with $r\geq 3$ for the Stirling numbers
$$\sum_{1\leq i_1<i_2<\dots<i_r\leq p-1} i_1i_2\dots i_r$$
or for their conjugates
$$\sum_{1\leq i_1<i_2<\dots<i_r\leq p-1} \frac{1}{i_1i_2\dots i_r},$$
the so-called multiple harmonic sums is both difficult and yet unknown in the general case. \\
Multiple harmonic sums have drawn much attention from mathematicians because they are part of a larger framework, namely they are special values of multiple zeta functions at positive integers. \\
In \cite{ZO}, Jianqiang Zhao shows that
$$\sum_{1\leq i<j\leq p-1}\frac{1}{ij}=\frac{\mh_{p-1}}{p}\mpt $$
but no systematic study gets performed by him nor by other authors. \\
In the congruence above, $\mh_{p-1}$ denotes the harmonic number of order $(p-1)$, sometimes denoted by $H_{p-1,1}$ within a context of generalized harmonic numbers. Wolstenholme's theorem asserts that $\mh_{p-1}$ is divisible by $p^2$. Sun's result (see point $(a)$ of Result $4$ in forthcoming section $\S\,1.2$) asserts that $H_{p-1,1}=p^2\mb_{p-3}\mpt$. By confronting our method and Zhao's result, we will find out $H_{p-1,1}$ modulo $p^4$.
\newtheorem*{Thm}{Theorem}
\begin{Thm}$\!\!\!0.0.$
$$H_{p-1,1}=-p^2(\mb_{2p-4}-2\mb_{p-3})\;\;\text{mod}\,p^4$$
\end{Thm}
\noindent Note, the modulus $p^3$ case is consistent with K\"ummer's congruences (see forthcoming Result $5$ of $\S\,1.2$).
Theorem $0.0$ is a special case of Sun's more general result.
In his ground breaking paper \cite{SU2}, Sun reveals without proof the generalized harmonic numbers $H_{p-1,k}$ modulo $p^4$. His result gets stated in his Remark $5.1.$

From our discussion earlier, the knowledge of $\as_2$ modulo $p^3$ provides the knowledge of the Stirling number $A_{p-3}$ modulo $p^3$. The statement appears in the following theorem, where $q_a$ denotes the Fermat quotient in base $a$.
\begin{Thm}$\!\!\!0.1.$
$$\left[\begin{array}{l}p\\3\end{array}\right]=(p\,B_{p-1}-p)\,p\,\mb_{p-3}-\frac{p^2}{2}\sum_{a=1}^{p-1}\frac{q_a^2}{a^2}\;\;\mpt$$
\end{Thm}
As we will see, a consequence of Theorem $0.1$ is the following, where $w_p$ denotes the Wilson quotient.
\begin{Thm}$\!\!\!0.2.$
$$\mc\mb(p-3):=\sum_{i=2}^{p-5}\mb_i\mb_{p-3-i}=2\,w_p\,\mb_{p-3}-\sum_{a=1}^{p-1}\frac{q_a^2}{a^2}\qquad\mpu$$
\end{Thm}
As part of our work here, we will determine $\as_4$.
\begin{Thm}$\!\!\!0.3.$ Using the notations of forthcoming Definition $1$, we have:
\begin{eqnarray*}
\as_4&=&p\,\mb_{p-5}+p^2\Big(\frac{1}{2}\mb_{p-3}^2-(\mb_{2(p-3)}-\mb_{p-5})_1\Big)\qquad\qquad\qquad\;\;\,\mpt\\
\left[\begin{array}{l}p\\5\end{array}\right]&=&-p\,\mb_{p-5}+p^2\Big(-\frac{1}{2}\mb_{p-3}^2+w_p\,\mb_{p-5}+\big(\mb_{2(p-3)}-\mb_{p-5}\big)_1\Big)\mpt\\
\mc\mb(p-5)&:=&\sum_{i=2}^{p-7}\mb_i\mb_{p-5-i}=-\mb_{p-3}^2+2\,w_p\mb_{p-5}+2(\mb_{2(p-3)}-\mb_{p-5})_1\;\mpu
\end{eqnarray*}
\end{Thm}
Outside of this study by Zhao for $\as_2$, part of the progress on the topic of multiple harmonic sums so far appears to be the existence of congruences that have become commonly referred to as "curious congruences on multiple harmonic sums", once one of the authors had qualified these congruences of "curious" and named them this way.
These congruences originate in \cite{ZO} when Zhao proves that for any prime $p\geq 5$,
$$\sum_{\begin{array}{l}i+j+k=p\\i,j,k>0\end{array}}\frac{1}{ijk}=-2\,B_{p-3}\mpu$$
In \cite{CH}, Chun-Gang Ji provides an independent proof using some combinatorial techniques.
Their result got generalized in \cite{ZC} to
$$\sum_{\begin{array}{l}l_1+l_2+\dots+l_n=p\\l_1,\dots,l_n>0\end{array}}\frac{1}{l_1l_2\dots l_n}=\begin{cases}
-(n-1)!\,B_{p-n} & \mpu\;\;\text{if $n$ is odd}\\
-\frac{n!\,n\,p}{2\,(n+1)}\,B_{p-n-1}&\mpd\;\;\text{if $n$ is even}\end{cases}$$
The most recent and elaborate versions go to higher powers of $p$ but only in the case of reciprocal products of $2$, $3$ or $4$ integers that are prime to $p$. On one hand in \cite{WC}, the authors show that for every positive integer $r$ and prime $p>2$,
$$\sum_{\begin{array}{l}i+j+k=p^r\\(i,j,k)\wedge p=1\end{array}}\frac{1}{ijk}=-2\,p^{r-1}\,B_{p-3}\;\;\text{mod}\;p^r$$
On the other hand in \cite{ZA}, Zhao considers the case of reciprocal products of $2$ or $4$ integers. He shows that for every positive integer
$r\geq \frac{n}{2}$ with $n\in\lbrace 2,4\rbrace$, and prime $p\geq 5$, we have
$$\sum_{\begin{array}{l}i_1+\dots+i_n=p^r\\(i_1,\dots,i_n)\wedge p=1\end{array}}\frac{1}{i_1i_2\dots i_n}=-\frac{n!}{n+1}\,p^r\,B_{p-n-1}\;\;\text{mod}\;p^{r+1}$$
Later on in \cite{WA2}, Liuquan Wang generalizes Jianqiang Zhao's result to sums involving reciprocal products of $6$ integers, stating that
for any prime $p\geq 11$ and any integer $r\geq 2$,
$$\sum_{\begin{array}{l}l_1+l_2+\dots+l_6=p^r\\(l_1,\dots,l_6)\wedge p=1\end{array}}\frac{1}{l_1l_2l_3l_4l_5l_6}=-\frac{5!}{18}\,p^{r-1}\,B_{p-3}^2\;\text{mod}\,p^r$$
The author had first studied the case of five variables in \cite{WA1} namely showing that
$$\sum_{\begin{array}{l}l_1+l_2+\dots+l_5=p^r\\(l_1,\dots,l_5)\wedge p=1\end{array}}\frac{1}{l_1l_2l_3l_4l_5}=-\frac{5!}{6}\,p^{r-1}\,B_{p-5}^2\;\text{mod}\,p^r$$
$$\begin{array}{l}\end{array}$$
Closing this digression on these "curious congruences", we will now state our results. Before that, we introduce some notations.
\newtheorem{Definition}{Definition}
\begin{Definition} (Hensel's $p$-adic expansion).
If $x$ is any $p$-adic integer, we will denote the coefficients of its expansion by the $(x)_i$'s as in
$$x=\sum_{i=0}^{+\infty}(x)_i\,p^i\in\mathbb{Z}_p$$
By times, we may also use a notation of \cite{LEV}, which consists of writing $x^{(i)}$ instead of $(x)_i$.
\end{Definition}
\begin{Definition} (K\"ummer prime and K\"ummer pair).
We say that $p$ is a K\"ummer prime if $\mb_{2(p-1)-2n}=\mb_{p-1-2n}\mpd$ for some integer $n$ with $2\leq 2n\leq p-3$. For such an $n$ when it exists, we will say that $(p,p-1-2n)$ is a K\"ummer pair.
\end{Definition}
\newtheorem{Remark}{Remark}
\begin{Remark}
In \cite{EM}, Ernvall and Mets\"ankyl\"a show that if $(p,p-1-2n)$ is an irregular pair, then $(p,p-1-2n)$ is not a K\"ummer pair.
\end{Remark}
\begin{Definition} (Truncated convolution of divided Bernoulli numbers). \\Let $n$ be an integer with $4\leq 2n\leq p-3$. We define:
$$\mt\mc\mb(p+1-2n, p-3):=\sum_{i=p+1-2n}^{p-3}\mb_i\mb_{2(p-1)-2n-i}$$
\end{Definition}
\begin{Definition} (Multiple harmonic sums). Define,
$$\as_k:=\sum_{1\leq i_1<\dots<i_k\leq p-1}\frac{1}{i_1\dots i_k}$$
\end{Definition}
\noindent Our main results get listed below.
\newtheorem{Theorem}{Theorem}
\begin{Theorem} \hfill\\
(i) Assume that $4\leq 2n\leq p-7$.
\begin{equation*}
\begin{split}
\frac{p^2}{2}\sum_{i=p+1-2n}^{p-3}\mb_i\mb_{2(p-1)-2n-i}=&-\frac{p^2}{2}\sum_{i=2}^{p-1-2n-2}\mb_i\mb_{p-1-2n-i}\\
&+p\,\big(\mb_{2(p-1)-2n}-\mb_{p-1-2n}\big)\\
&+p^2\bigg((p\,B_{p-1})_1-1\bigg)\,\mb_{p-1-2n}\mpt
\end{split}
\end{equation*}
(ii) Case $2n=p-3$.
\begin{equation*}
p^2\sum_{i=4}^{p-3}\mb_i\mb_{p+1-i}=2p\,\mb_{p+1}+p^2\mb_2+2p\mb_2(p\,B_{p-1})\mpt
\end{equation*}
(iii) Case $2n=p-5$.
\begin{equation*}
p^2\sum_{i=6}^{p-3}\mb_i\mb_{p+3-i}=\frac{7}{720}\,p^2+2p\,\mb_{p+3}+2p\,\mb_4(p\,B_{p-1})\mpt
\end{equation*}
\end{Theorem}
\noindent We note that by Glaisher's result, stated as Result $1$ in $\S\,1.2$, the Wilson quotient $w_p:=\frac{(p-1)!+1}{p}\mpu$ is equal to $(p\,B_{p-1})_1-1$. In particular, we see that if $p$ is a Wilson prime and $(p,p-1-2n)$ is a K\"ummer pair, then the truncated convolution
$\mt\mc\mb(p+1-2n,p-3)$ is congruent modulo $p\,\mathbb{Z}_p$ to the opposite of the full convolution of order $p-1-2n$ when $4\leq 2n\leq p-7$. \\
In \cite{CDP}, the authors showed that there are only three Wilson primes less than $5\times 10^8$, namely $5$, $13$ and $563$. It is unknown whether there exists any newly defined K\"ummer prime and even more unknown whether there exist any Wilson primes that are also K\"ummer primes.\\

Our next theorem deals with multiple harmonic sums modulo $p^3$.
\begin{Theorem}\hfill\\
(I) Multiple harmonic sums \textbf{modulo} $\mathbf{p^3}$.
$$\as_k=\begin{cases}\frac{k+1}{2}\,p^2\,\mb_{p-2-k}&\text{if $k$ is odd and $k\leq p-4$}\\
&\\
\frac{p}{2}-p^2+\frac{1}{2}(p\,B_{p-1})_1p^2&\text{if $k=p-2$}\\
&\\
\,p\,\big(2\,\mb_{p-1-k}-\mb_{2(p-1)-k}\big)+\frac{p^2}{2}\mt\mc\mb(p+1-k,p-3)&\text{if $k$ is even and $k\leq p-5$}\\
&\\
\frac{p}{12}-\frac{11p^2}{24}+\frac{p^2}{12}(p\,B_{p-1})_1&\text{if $k=p-3$}\\
&\\
3(p\,\mb_{p-1}-1)-p\,\mb_{2(p-1)}-\frac{1}{2}p^2\,\mb_{p-1}^2&\text{if $k=p-1$}
\end{cases}$$
(II) Let $w_p$ denote the Wilson quotient, namely $w_p=(p\,B_{p-1})_1-1$.\\Suppose $k\in\lbrace 1,\,\dots,p-2\rbrace$. The following congruence holds.
\begin{equation*}
\as_{k}=-A_{p-1-k}+(-1)^k\,w_p\,p^2\,\mb_{p-1-k}\mpt
\end{equation*}
(III) Let $j\in\lbrace 1,\dots,p-2\rbrace$. The $\as_j$'s are all divisible by $p$ and
\begin{equation*}
A_{p-1-j}=(p\,B_{p-1}-p)\as_j\mpt
\end{equation*}
\end{Theorem}

\subsection{Notations and some useful results}
Throughout the paper it will be useful to know the following generalization of Wilson's theorem, established by Glaisher in $1900$ \cite{GL}. See also \cite{LEV} for an independent proof of that result.
\newtheorem{Result}{Result}
\begin{Result}(Wilson's theorem modulo $p^2$, Glaisher \cite{GL}, $1900$)
$$(p-1)!=p\,B_{p-1}-p\mpd$$
\end{Result}
Some sums of importance throughout the paper are the sums of powers of the first $(p-1)$ integers and the generalized harmonic numbers $H_{p-1,k}$ which for convenience we sometimes simply denote by $H_k$:
$$S_k=\sum_{a=1}^{p-1}a^k\qquad\text{and}\qquad H_k=\sum_{a=1}^{p-1}\frac{1}{a^k}$$
We are especially interested in congruences concerning these sums modulo powers of $p$.
A result of \cite{SU2} reads
\begin{Result} (Sun \cite{SU2}, $2000$) Let $k$ be an integer with $k\geq 2$. Then, we have:
$$S_k=p\,B_k+\frac{p^2}{2}\,k\,B_{k-1}+\frac{p^3}{6}\,k(k-1)\,B_{k-2}\;\;\text{mod}\,p^3$$
\end{Result}
\noindent An extensively used result on the Bernoulli numbers dating from the $19th$ century and due to Von Staudt \cite{VS} and independently Clausen \cite{CL}, implies that the denominators of the Bernoulli number $B_{2m}$ consists of a product of primes $q$ with multiplicities one, such that $q-1$ divides $2m$. This is particularly useful to know when dealing with congruences modulo powers of $p$. \\
Regarding the harmonic numbers, a classical result is the Wolstenholme's theorem \cite{WO} dating from $1862$ which asserts that
\begin{center}$H_{p-1,1}=0\mpd$ and $H_{p-1,2}=0\mpu$\end{center}
The result got generalized by Bayat in $1997$ in \cite{BA}. Bayat deals with other $k$'s as well with $k\geq 3$. His result is stated below.
\begin{Result}(Bayat's generalization of Wolstenholme's theorem, 1997 \cite{BA}).\\
Let $m$ be a positive integer and let $p$ be a prime number with $p\geq m+3$. Then,
$$\sum_{k=1}^{p-1}\frac{1}{k^m}\;\equiv\;\begin{cases}0\;\;\text{mod}\;p\;\;\;\;\text{if $m$ is even}\\
0\;\;\text{mod}\;p^2\;\;\text{if $m$ is odd}
\end{cases}$$
\end{Result}
We also recall below one of the main results of \cite{SU2}.
\begin{Result} (Sun \cite{SU2}, $2000$). Let $p$ be a prime greater than $3$. Then,\\
(a) If $k\in\lbrace 1,2,\dots,p-4\rbrace$ then,
$$\sum_{a=1}^{p-1}\frac{1}{a^k}=\begin{cases}\frac{k(k+1)}{2}\frac{B_{p-2-k}}{p-2-k}p^2\qquad\qquad\;\mpt&\text{if $k$ is odd}\\
k\bigg(\frac{B_{2p-2-k}}{2p-2-k}-2\frac{B_{p-1-k}}{p-1-k}\bigg)p\;\mpt&\text{if $k$ is even}
\end{cases}$$
(b) $$\sum_{a=1}^{p-1}\frac{1}{a^{p-3}}=\bigg(\frac{1}{2}-3B_{p+1}\bigg)p-\frac{4}{3}p^2\mpt$$
(c) $$\sum_{a=1}^{p-1}\frac{1}{a^{p-2}}=-\big(2+p\,B_{p-1}\big)p+\frac{5}{2}p^2\mpt$$
(d) $$\sum_{a=1}^{p-1}\frac{1}{a^{p-1}}=p\,B_{2p-2}-3p\,B_{p-1}+3(p-1)\mpt$$
\end{Result}

We will also use the well-known K\"ummer congruences from \cite{KU} which we recall below. We also state the generalization by Sun in \cite{SU2} to the modulus $p^2$.
\begin{Result}(K\"ummer's congruences, $1850$ \cite{KU}, followed by Sun's generalization, $2000$ \cite{SU2})\\
Let $p$ be an odd prime and let $b>0$ be an even integer such that $p-1\not|\,b$. Then, we have for all nonnegative integer $k$,
$$\frac{B_{k(p-1)+b}}{k(p-1)+b}\equiv\frac{B_b}{b}\;\;\text{mod}\,p$$
$$\frac{B_{k(p-1)+b}}{k(p-1)+b}\equiv k\frac{B_{p-1+b}}{p-1+b}-(k-1)(1-p^{b-1})\frac{B_b}{b}\;\;\text{mod}\,p^2$$
\end{Result}
When $k=1$, Sun's generalization is a trivial congruence. While studying the irregular primes, Ernvall and Mets\"ankyl\"a have shown the following congruence, where $q_a$ denotes the Fermat quotient in base $a$.
\begin{Result} (Ernvall and Mets\"ankyl\"a \cite{EM}, $1991$)
Let $n$ be an even integer with $4\leq n\leq p-3$. Then, we have:
\begin{equation*}
\mb_{p-1+n}=\mb_{n}-\frac{p}{2}\,\sum_{a=1}^{p-1}q_a^2\,a^n\;\;\mpd
\end{equation*}
\end{Result}

A useful version of Result $4$ modulo $p^2$ for the even indices is originally due to Glaisher and can be derived from his results stated as Result $2$ and Theorem $4$ in \cite{LEV}. In the case when $2k\leq p-5$, it can also be derived from Sun's version in the second row of $(a)$ in Result $4$ by using the K\"ummer congruences above with integers $1$ and $p-1-k$ (recall that $k$ is even). We state this useful version as Result $7$ below.
\begin{Result}Let $k$ be an integer with $1\leq k\leq \frac{p-3}{2}$. Then, we have
$$H_{2k}=\frac{2k}{2k+1}p\,B_{p-1-2k}\mpd$$
\end{Result}
We will sometimes denote the unsigned Stirling numbers of the first kind $\left[\begin{array}{l}p\\s\end{array}\right]$ by $A_{p-s}$, using Glaisher's notation of \cite{GL}. We now recall a result of \cite{LEV}.
\begin{Result}(Stirling numbers modulo $p^3$, cf \cite{LEV}, preprint $2019$). \\Let $k$ be an integer with $1\leq k\leq \frac{p-1}{2}$. We have,
\begin{eqnarray*}
(k\leq\frac{p-3}{2})\;A_{2k+1}\negthickspace\negthickspace&=&\negthickspace\!\!\frac{p^2}{2}\frac{2k+1}{2k}\,B_{2k}\qquad\qquad\qquad\qquad\qquad\!\text{mod}\,p^3\\
&&\notag\\
(k\neq 1)\;A_{2k}\negthickspace\negthickspace&=&\negthickspace\!\!-\frac{1}{2k}\Bigg(p\,B_{2k}-p^2\,\sum_{r=1}^{ k-1}\frac{B_{2r}B_{2k-2r}}{2r}\Bigg)\;\;\;\text{mod}\,p^3\\
&&\notag\\
\text{and}\;A_1\negthickspace\negthickspace&=&\negthickspace\!\!\frac{p(p-1)}{2}\qquad\qquad\qquad\qquad\qquad\qquad\text{mod}\,p^3\\
&&\notag\\
(p\geq 5)\;\text{and}\;A_2\negthickspace\negthickspace&=&\negthickspace\!\!\frac{1}{2}\bigg(-\frac{p}{6}+\frac{3\,p^2}{4}\bigg)\qquad\qquad\qquad\qquad\;\text{mod}\,p^3
\end{eqnarray*}
\end{Result}
\noindent One of the formulas above involves a convolution of divided Bernoulli numbers with Bernoulli numbers. Convolutions involving Bernoulli numbers have drawn the interest of mathematicians over the centuries. Already Euler had found an identity.
\begin{Result}(Euler's identity, see e.g. \cite{SD})
$$\forall\;n\geq 1,\;\sum_{j=0}^n\binom{n}{j}B_jB_{n-j}=-n\,B_{n-1}-(n-1)B_n$$
\end{Result}
Japanese mathematician Hiroo Miki comes up with an identity in $1978$ which involves both a binomial convolution and an ordinary convolution of divided Bernoulli numbers, see \cite{MI}. Denoting
$$\mathcal{B}_n=\frac{B_n}{n},\;\;\text{and}\;\;\mathcal{H}_n=1+\frac{1}{2}+\dots+\frac{1}{n},$$
his identity reads:
\begin{Result}(Miki's identity \cite{MI}, $1978$)
$$\forall n>2,\,\sum_{i=2}^{n-2}\mathcal{B}_i\mathcal{B}_{n-i}=\sum_{i=2}^{n-2}\binom{n}{i}\mathcal{B}_i\mathcal{B}_{n-i}+2\mathcal{H}_n\mathcal{B}_n$$
\end{Result}
Miki shows that both sides of the identity are congruent modulo $p$ for sufficiently large $p$, which implies that they are equal. A more elementary proof of Miki's identity was given more recently by Ira Gessel \cite{GE}, using two different expressions for the Stirling numbers of the second kind. \\
Another version of Miki's identity, whose proof is inspired from quantum field theory, is due to Gerald V. Dunne and Christian Schubert and is the following.
\begin{Result} (Modified form of Miki's identity by Dunne and Schubert \cite{DS}, $2004$)
$$\sum_{k=1}^{n-1}\frac{B_{2k}B_{2n-2k}}{2k(2n-2k)}=\frac{1}{n}\sum_{k=1}^{n-1}\frac{B_{2k}B_{2n-2k}}{2k}\binom{2n}{2k}+\frac{B_{2n}}{n}H_{2n}$$
\end{Result}
\noindent Other Miki type identities, for instance concerning non divided Bernoulli numbers, get listed in \cite{DI}.

Along the paper, we will use Vandermonde's binomial convolution which we recall below.
\begin{Result}(Chu-Vandermonde's convolution, Chu Shi-Chieh $1303$ and Alexandre-Th\'eophile Vandermonde $1772$).\\
Let $m$, $n$ and $r$ be non-negative integers. Then,
$$
\binom{m+n}{r}=\sum_{k=0}^r\binom{m}{k}\binom{n}{r-k}
$$
\end{Result}

Last, we will make use of the polynomial
$$f=X^{p-1}+(p-1)!\in\mathbb{Z}_p[X]$$
For an integer $k$ such that $1\leq k\leq p-1$, we will define, using the same notations as in \cite{LEV}, $\delta_0(k)$ and $\delta_1(k)$ so that: $$k^{p-1}=1+p\,\delta_0(k)+p^2\delta_1(k)\mpt$$
Thus, $\delta_0(k)$ is the residue of the Fermat quotient $q_p(k)$ modulo $p$. \\
It is shown in \cite{LEV} that the polynomial $f$ factors as:
\begin{Result}(Factorization of $X^{p-1}+(p-1)!$ in $\mathbb{Z}_p[X]$, cf \cite{LEV}, $2019$)
$$f=(X-1-pt_1)(X-2-pt_2)\dots (X-(p-1)-pt_{p-1}),$$
with $t_i=\sum_{s=0}^{\infty}t_i^{(s)}\in\mathbb{Z}_p$ and
\begin{eqnarray*}
pt_i^{(0)}&=&i(1+(p-1)!+p\,\di)\mpd\\
t_i^{(1)}&=&i\bigg(\di+\dui+\big(\sum_{k=1}^{p-1}\dk\big)^2+(1+\di)\sum_{k=1}^{p-1}\dk\bigg)\;\text{mod}\,p
\end{eqnarray*}
\end{Result}

\section{The Stirling numbers modulo $p^3$}
\subsection{Stirling numbers with even indices}
We first compute the Stirling numbers with even indices modulo $p^3$. Due to the past works of Glaisher and Sun, our method will show relevant for the odd indices only, but we start with the even indices anyway in order to get familiar with it.\\
By definition, $\pdn$ is the unsigned coefficient of $X^{2n-1}$ in
$$(X-1)(X-2)\dots (X-(p-1))$$
Because $\left[\begin{array}{l}\;\;\,p\\p-1\end{array}\right]$ can be computed directly and is simply $\frac{p(p-1)}{2}$ modulo $p^3$,
from now on we will focus on the even indices $2n$ such that $2\leq 2n\leq p-3$. \\
By \cite{LEV}, we have the factorization in $\mathbb{Z}_p[X]$,
$$X^{p-1}-1=(X-1-pt_1)(X-2-pt_2)\dots (X-(p-1)-pt_{p-1}),$$
with the first two coefficients of the $p$-adic expansion of the $t_i$'s provided in \cite{LEV} and recalled here further below. In what follows, we will denote this polynomial by $f(X)$ or simply $f$. We look at the coefficient in $X^{2n-1}$ in both factored and expanded form of $f$ modulo $p^3$.
We have, where the range of integers lies in $\lbrace1,\dots,p-1\rbrace$,
\hspace{-1in}\begin{equation}\begin{split}\pdn=
&-\frac{1}{(2n-1)!}\sum_{m_0}p\,t_{m_0}^{(0)}(p-1)!^{m_0}\sum_{m_1\neq m_0}\frac{1}{m_1}\sum_{m_2\neq m_1,m_0}\frac{1}{m_2}\dots\sum_{m_{2n-1}\neq m_0,\dots,m_{2n-2}}\frac{1}{m_{2n-1}}\\
&-\frac{1}{(2n-1)!}\sum_{m_0}p^2\,t_{m_0}^{(1)}(p-1)!^{m_0}\sum_{m_1\neq m_0}\frac{1}{m_1}\sum_{m_2\neq m_1,m_0}\frac{1}{m_2}\dots\sum_{m_{2n-1}\neq m_0,\dots,m_{2n-2}}\frac{1}{m_{2n-1}}\\
&-\frac{(p-1)!}{2(2n-1)!}\sum_{m_0}\frac{pt_{m_0}^{(0)}}{m_0}\sum_{m_1\neq m_0}\frac{pt_{m_1}^{(0)}}{m_1}\sum_{m_2\neq m_0,m_1}\frac{1}{m_2}\dots
\sum_{m_{2n}\neq m_0,m_1,\dots, m_{2n-1}}\frac{1}{m_{2n}}\mpt
\end{split}\end{equation}
We will denote the first (resp second, resp third) row by $\mathcal{S}_1$ (resp $\mathcal{S}_2$, resp $\mathcal{S}_3$). We will deal with each row independently.\\ Thus, we have:
\begin{equation}
\pdn=\mathcal{S}_1+\mathcal{S}_2+\mathcal{S}_3\mpt
\end{equation}
We proceed the sums from right to left by first summing over all the $(p-1)$ terms and then subtracting the exceeding terms. We will show the following lemma.
\newtheorem{Lemma}{Lemma}
\begin{Lemma}We have,
\begin{eqnarray*}(i)\;\;\mathcal{S}_1&=&H_{2n-1}-S_{p-2n}+p^2\sum_{m_0=1}^{p-1}\frac{\dmzu}{m_0^{2n-1}}\mpt\\
(ii)\;\;\mathcal{S}_2&=&-p^2\sum_{m_0=1}^{p-1}\frac{\dmzu}{m_0^{2n-1}}\qquad\qquad\qquad\;\;\;\mpt\\
(iii)\;\;\mathcal{S}_3&=&n\big(S_{2p-2n-1}-2S_{p-2n}+H_{2n-1}\big)\;\;\,\mpt
\end{eqnarray*}
\end{Lemma}
\newtheorem{Corollary}{Corollary}
\begin{Corollary}We have,
$$\pdn=(n+1)\,H_{2n-1}+n\,S_{2p-2n-1}-(2n+1)\,S_{p-2n}\mpt$$
\end{Corollary}
\textsc{Proof of Lemma $1$}.
We first deal with $\mathcal{S}_1$. By Bayat's result and since we work modulo $p^3$, there can be at most one "full" sum when we consider the last $2n-1$ sums. Moreover, we will distinguish between "exactly one full sum" and "no full sum". In the latter case, we get the contribution:
\begin{equation}(p-1)!\sum_{m_0=1}^{p-1}\frac{p\,t_{m_0}^{(0)}}{m_0^{2n}}\end{equation}
We now examine "exactly one full sum". A general form for the contribution is:
\begin{equation}
-\frac{(p-1)!}{(2n-1)!}\sum_{s=1}^{n-1}C_{2s}\sum_{m=1}^{p-1}\frac{1}{m^{2s}}\sum_{m_0=1}^{p-1}\frac{p\,t_{m_0}^{(0)}}{m_0^{2n-2s}}
\end{equation}
for some integer coefficients $C_{2s}$'s. There is actually no need to evaluate these integer coefficients as we can show directly that $(4)$ is congruent to zero modulo $p^3\mathbb{Z}_p$. Indeed, we have by \cite{LEV},
$$pt_{m_0}^{(0)}=m_0(1+(p-1)!+p\,\delta_0(m_0))\mpd$$ with
$$m_0^{p-1}=1+p\delta_0(m_0)+p^2\dmzu\mpt$$
Then, replacing in $(4)$ yields:
$$-\frac{(p-1)!}{(2n-1)!}\sum_{s=1}^{n-1}C_{2s}H_{2s}\sum_{m_0=1}^{p-1}\frac{1+(p-1)!+p\,\delta_0(m_0)}{m_0^{2n-2s-1}}$$
Next, by Bayat's theorem, we know that
$$H_{2s}=0\mpu\;\text{and}\;H_{2n-2s-1}=0\mpd$$
Also, by Congruence $(5.1)$ of \cite{SU2}, we have
$$S_{p-2n+2s}=0\mpd$$
It follows that $(4)$ is congruent to zero modulo $p^3$. \\
It remains to evaluate $(3)$ modulo $p^3$. Modulo $p^3$, expression $(3)$ reduces to
$$(p-1)!\bigg(S_{p-2n}-H_{2n-1}-p^2\sum_{m_0=1}^{p-1}\frac{\dmzu}{m_0^{2n-1}}\bigg),$$
which in turn reduces to $$H_{2n-1}-S_{p-2n}+p^2\sum_{m_0=1}^{p-1}\frac{\dmzu}{m_0^{2n-1}},$$ after application of Wilson's theorem since
$$S_{p-2n}=0\mpd\;\text{and}\;H_{2n-1}=0\mpd$$
This settles point $(i)$ of Lemma $1$. Next, we note that $\mathcal{S}_2$ is obtained from $\mathcal{S}_1$ by replacing $p\,t_{m_0}^{(0)}$ by $p^2\,t_{m_1}^{(0)}$.
Then, there is no full sum to be kept during the evaluation. \\We thus get:
$$\mathcal{S}_2=-\sum_{m_0=1}^{p-1}\frac{p^2t_{m_0}^{(1)}}{m_0^{2n}}\mpt$$
Moreover, we have by Lemma $3$ of \cite{LEV},
$$t_{m_0}^{(1)}=m_0\bigg(\delta_0(m_0)+\delta_1(m_0)+\big(\sum_{i=1}^{p-1}\di\big)^2+(1+\dmz)\sum_{i=1}^{p-1}\di\bigg)\mpu$$
After inspection, there is only one term contributing to $\mathcal{S}_2$ modulo $p^3$. Namely we get,
$$\mathcal{S}_2=-p^2\sum_{m_0=1}^{p-1}\frac{\dmzu}{m_0^{2n-1}}\mpt,$$
which constitutes point $(ii)$ of Lemma $1$.\\
It remains to deal with $\mathcal{S}_3$. No sum to the right when excluding the first two sums may contribute fully. Also, if the second sum to the left contributes fully, then we get a term in
$$\sum_{m_0}\frac{pt_{m_0}^{(0)}}{m_0^r}\sum_{m_1}\frac{pt_{m_1}^{(0)}}{m_1^s}$$
with $r+s=2n+1$. In particular, we see that $r$ and $s$ have distinct parity, hence after investigation, the term is congruent to zero modulo $p^3$.
Consequently, we simply obtain
$$\mathcal{S}_3=-n(p-1)!\sum_{m_0}\frac{p^2(t_{m_0}^{(0)})^2}{m_0^{2n+1}}\mpt$$
It then reduces to
$$\mathcal{S}_3=-n(p-1)!\sum_{m_0}\frac{p^2\,\dmz^2}{m_0^{2n-1}}\mpt$$
And so, we get
$$\mathcal{S}_3=-n(p-1)!\big(S_{2p-2n-1}-2S_{p-2n}+H_{2n-1}\big)\mpt$$
All the indices inside the parenthesis are odd, hence $(p-1)!$ must be taken modulo $p$. We thus obtain point $(iii)$ of Lemma $1$.\\\\
We immediately derive Corollary $1$. \\\\
Further, by Congruence $(5.1)$ of \cite{SU2} stated as Result $2$ in the introduction, we have when $2\leq 2n\leq p-3$,
\begin{eqnarray*}S_{2p-2n-1}&=&-(2n+1)\frac{p^2}{2}\,B_{2p-2n-2}\;\mpt\\
S_{p-2n}&=&-np^2\,B_{p-2n-1}\;\qquad\;\;\;\;\;\,\mpt
\end{eqnarray*}
Also, by Theorem $5.1$ of \cite{SU2} stated as Result $4$ in the introduction, we have when $2\leq 2n\leq p-3$,
$$H_{2n-1}=
-\frac{n(2n-1)}{2n+1}\,p^2B_{p-2n-1}\;\;\mpt$$
Then, Corollary $1$ rewrites as:
\begin{equation}\pdn=n\,\frac{2n^2+3n+2}{2n+1}p^2\,B_{p-2n-1}-\frac{n(2n+1)}{2}p^2\,B_{2p-2n-2}\mpt\end{equation}
Moreover, by K\"ummer's congruences, we have since $p-1$ does not divide $2n$:
$$\frac{B_{p-1+p-2n-1}}{p-1+p-2n-1}=\frac{B_{p-2n-1}}{p-2n-1}\;\text{mod}\,p$$
Hence, it comes:
\begin{equation}B_{2p-2n-2}=\frac{2(n+1)}{2n+1}\,B_{p-2n-1}\;\text{mod}\,p\end{equation}
Gathering $(5)$ and $(6)$, it follows that:
$$\pdn=\frac{n}{2n+1}p^2\,B_{p-2n-1}\mpt$$
This is precisely Glaisher's result from \cite{GL} which is stated as Result $4$ in \cite{LEV} using Glaisher's notations and which is also listed in Result $8$ of the introduction.

\subsection{Stirling numbers with odd indices}
We now deal with the case of interest, namely compute the Stirling numbers modulo $p^3$ with odd indices. We adapt the work we did with the even indices. This case is more complicated due to the fact that the generalized harmonic numbers with even indices are only congruent to zero modulo $p$ instead of being congruent to zero modulo $p^2$ like their homologs with the odd indices (cf Bayat's generalization of Wolstenholme's theorem). Likewise, the sums of powers with even powers are no longer divisible by $p^2$, unlike their homologs with the odd powers.\\
First, we will get some general expressions involving some convolutions and inside them some yet unknown integer coefficients. Second, we will do the combinatorial analysis of these coefficients. Third, we will resolve the convolutions. \\

We start below with the first point.
\subsubsection{General calculation}
By looking at the coefficient of $X^{2n}$ in the polynomial $f$ in both factored and expanded forms, we get modulo $p^3$:
\hspace{-1.5in}\begin{equation}\begin{split}0=\left[\begin{array}{l}\;p\;\\2n+1\end{array}\right]+
&\frac{(p-1)!}{(2n)!}\sum_{m_0}\frac{p\,t_{m_0}^{(0)}}{m_0}\sum_{m_1\neq m_0}\frac{1}{m_1}\sum_{m_2\neq m_1,m_0}\frac{1}{m_2}\dots\sum_{m_{2n}\neq m_0,\dots,m_{2n-1}}\frac{1}{m_{2n}}\\
&\negthickspace\negthickspace\negthickspace\negthickspace+\frac{(p-1)!}{(2n)!}\sum_{m_0}\frac{p^2\,t_{m_0}^{(1)}}{m_0}\sum_{m_1\neq m_0}\frac{1}{m_1}\sum_{m_2\neq m_1,m_0}\frac{1}{m_2}\dots\sum_{m_{2n}\neq m_0,\dots,m_{2n-1}}\frac{1}{m_{2n}}\\
&\negthickspace\negthickspace\negthickspace\negthickspace\negthickspace\negthickspace\negthickspace\negthickspace\negthickspace\negthickspace\negthickspace\negthickspace
\negthickspace\negthickspace\negthickspace\negthickspace+\frac{(p-1)!}{2(2n)!}\sum_{m_0}\frac{pt_{m_0}^{(0)}}{m_0}\sum_{m_1\neq m_0}\frac{pt_{m_1}^{(0)}}{m_1}\sum_{m_2\neq m_0,m_1}\frac{1}{m_2}\dots
\sum_{m_{2n+1}\neq m_0,m_1,\dots, m_{2n}}\frac{1}{m_{2n+1}}\\
&\\
&\qquad\qquad\qquad\qquad\qquad\qquad\qquad\qquad\qquad\qquad\qquad\qquad\mpt
\end{split}\end{equation}
Like before, the range of integers for the respective sums lies in $\lbrace 1,\dots,p-1\rbrace$. Also, we will denote the sum on the first (resp second, resp third) row by $\mathcal{S}^{'}_1$ (resp $\mathcal{S}^{'}_2$, resp $\mathcal{S}^{'}_3$), so that
$$\si=-\ms1-\ms2-\ms3\mpt$$ We show the following lemma.
\begin{Lemma}We have,\\
$(i)\;\;$ \begin{equation*}\begin{split}\ms1=(2p-1-2p\,B_{p-1})H_{2n}&+(pB_{p-1}-p)S_{p-1-2n}+\sum_{m_0}\frac{p^2\dmzu}{m_0^{2n}}\\
&+\frac{1}{(2n)!}\sum_{s=1}^{n-1} D_sH_{2s}(S_{p-1-2n+2s}-H_{2n-2s})\mpt\end{split}\end{equation*}
$(ii)\;\,$\begin{equation*}\ms2=(pB_{p-1}+1)(H_{2n}-S_{p-1-2n})-\sum_{m_0}\frac{p^2\dmzu}{m_0^{2n}}\mpt\end{equation*}
$(iii)\;\,$\begin{equation*}\begin{split}\ms3=&\frac{2n+1}{2}\bigg((2p-1-2p\,B_{p-1})H_{2n}+S_{2p-2-2n}+2(p\,B_{p-1}-p)S_{p-1-2n}\bigg)\\
&-\frac{1}{2(2n)!}\sum_{k=1}^{n-1} G_k(S_{p-1-2k}-H_{2k})(S_{p-1-2n+2k}-H_{2n-2k})\mpt\end{split}\end{equation*}
\end{Lemma}
\begin{Corollary}
\begin{equation*}\begin{split}\si=&\;\Big((2n+2)\,p\,B_{p-1}-(2n+3)p+\frac{2n+1}{2}\Big)H_{2n}\\&
\Big(-(2n+1)\,p\,B_{p-1}+(2n+2)p+1\Big)\,S_{p-1-2n}\\
&-\frac{2n+1}{2}\,S_{2p-2n-2}\\
&-\frac{1}{(2n)!}\sum_{s=1}^{n-1} D_sH_{2s}(S_{p-1-2n+2s}-H_{2n-2s})\\
&+\frac{1}{2(2n)!}\sum_{k=1}^{n-1} G_k(S_{p-1-2k}-H_{2k})(S_{p-1-2n+2k}-H_{2n-2k})\mpt\end{split}\end{equation*}
with the $D_s$'s and the $G_k$'s some integer coefficients to determine.
\end{Corollary}
\textsc{Proof of Lemma $2$}. We proceed like in the even case. There can be at most one full sum when processing $\ms1$.
\begin{enumerate}[(a)]
\item
If there is no full sum, we obtain the contribution modulo $p^3$,
\begin{equation}(p-1)!\sum_{m_0=1}^{p-1}\frac{p\,t_{m_0}^{(0)}}{m_0^{2n+1}}\end{equation}
Modulo $p^3$, the expression $(8)$ is congruent to
\begin{equation}
-(1+(p-1)!)\sum_{m_0}\frac{1}{m_0^{2n}}+(p-1)!\sum_{m_0}\frac{p\dmz}{m_0^{2n}},
\end{equation}
which in turn is congruent modulo $p^3$ to
\begin{equation}
(2p-1-2p\,B_{p-1})H_{2n}+(p\,B_{p-1}-p)\,S_{p-1-2n}+\sum_{m_0}\frac{p^2\dmzu}{m_0^{2n}},
\end{equation}
where we used Result $1$ from the introduction.
\item If there is exactly one full sum, the general form for the contribution modulo $p^3$ is
\begin{equation}
-\frac{(p-1)!}{(2n)!}\sum_{s=1}^n D_s\sum_{m}\frac{1}{m^{2s}}\sum_{m_0}\frac{pt_{m_0}^{(0)}}{m_0^{2n+1-2s}},
\end{equation}
with some integer coefficients $D_s$'s to be determined. \\
The term arising from $2s=2n$ yields a null contribution. Thus,
Expression $(11)$ is in turn congruent modulo $p^3$ to
\begin{equation}
\frac{1}{(2n)!}\sum_{s=1}^{n-1}D_s\,H_{2s}\bigg(S_{p-1-2n+2s}-H_{2n-2s}\bigg)
\end{equation}
\end{enumerate}
Point $(i)$ of Lemma $2$ is then obtained by gathering the results from points $(a)$ and $(b)$ above. \\

We now process $\ms2$. Using the terminology from before, there is no full sum to be kept. We must evaluate modulo $p^3$,
\begin{equation}
-\sum_{m_0=1}^{p-1}\frac{p^2\,t_{m_0}^{(1)}}{m_0^{2n+1}}
\end{equation}
Recall from \cite{LEV} that
$$(p-1)!=-1+p\sum_{k=1}^{p-1}\dk\mpd$$
This is Theorem $1$ of \cite{LEV}. It is obtained from using some $p$-adic techniques such as Hensel's lifting algorithm for lifting the $p$-adic integer roots of the polynomial $f$ to the first power of $p$.
Using this fact together with Result $13$, we obtain point $(ii)$ of Lemma $2$. \\

It remains to tackle $\ms3$. Again, there is no full sum to be kept concerning the last $2n$ sums since we work modulo $p^3$. \\\\
Now, there are two cases:
\begin{enumerate}[(c)]
\item Either we take the first two full sums to the left and we thus have a term in
$$\frac{(p-1)!}{2(2n)!}\sum_{m_0}\frac{pt_{m_0}^{(0)}}{m_0^r}\sum_{m_1}\frac{pt_{m_1}^{(0)}}{m_1^s},$$
with $r+s=2n+2$. We see that either $r$ and $s$ are both even or $r$ and $s$ are both odd. In the first case we get zero modulo $p^3$. Only the second case yields a non-zero contribution. More precisely, we get
\begin{equation}
\frac{(p-1)!}{2(2n)!}\sum_{k=0}^n G_k\sum_{m_0}\frac{pt_{m_0}^{(0)}}{m_0^{2k+1}}\sum_{m_1}\frac{pt_{m_1}^{(0)}}{m_1^{2n+2-2k-1}},
\end{equation}
for some adequate integer coefficients $G_k$'s.
Then, by using Result $13$, we get
\begin{equation}
\frac{(p-1)!}{2(2n)!}\sum_{k=1}^{n-1} G_k\sum_{m_0}\frac{1+(p-1)!+p\dmz}{m_0^{2k}}\sum_{m_1}\frac{1+(p-1)!+p\delta_0(m_1)}{m_1^{2n-2k}}
\end{equation}
From there, the only contribution that is left modulo $p^3$ is
\begin{equation}
\frac{(p-1)!}{2(2n)!}\sum_{k=1}^{n-1} G_k\sum_{m_0}\frac{p\dmz}{m_0^{2k}}\sum_{m_1}\frac{p\delta_0(m_1)}{m_1^{2n-2k}}
\end{equation}
And a new reduction modulo $p^3$ yields:
\begin{equation}
-\frac{1}{2(2n)!}\sum_{k=1}^{n-1} G_k(S_{p-1-2k}-H_{2k})(S_{p-1-2n+2k}-H_{2n-2k})
\end{equation}
\end{enumerate}
\begin{enumerate}[(d)]
\item Or we take no full sum except the very first one. We get modulo $p^3$,
\begin{equation}
\frac{2n+1}{2}\sum_{m_0}\frac{p^2\big(t_{m_0}^{(0)}\big)^2}{m_0^{2n+2}}
\end{equation}
It yields after reduction modulo $p^3$:
\begin{equation}
\frac{2n+1}{2}\bigg((2p-1-2p\,B_{p-1})H_{2n}+S_{2p-2-2n}+2(p\,B_{p-1}-p)S_{p-1-2n}\bigg)
\end{equation}
\end{enumerate}
The results obtained in $(c)$ and $(d)$ are gathered in Point $(iii)$ of Lemma $2$. \\
The next step is to calculate the integer coefficients $D_s$'s and $G_k$'s.

\subsubsection{Combinatorial analysis of the coefficients $D_s$'s and $G_k$'s}

\begin{Lemma}\hfill\\

\noindent (i) Let $s$ be an integer with $1\leq s\leq n-1$. We have,
\begin{equation}
D_s=\frac{(2n)!}{2s}
\end{equation}
(ii) Let $k$ be an integer with $1\leq k\leq n-1$. We have,
\begin{equation}
G_k=(2n)!
\end{equation}
\end{Lemma}
\textsc{Proof of Lemma $3$}. \hfill\\

The integer coefficient $D_s$ counts the number of ways to fix one full sum with even power $2s$, e.g., $\sum_{m_i}\frac{1}{m_i^{2s}}$, such that the withdrawal process like explained before on all the other indices in
$$\sum_{m_1}\frac{pt_{m_1}^{(0)}}{m_1}\sum_{m_2\neq m_1}\frac{1}{m_2}\dots\sum_{m_i}\frac{1}{m_i}\dots\sum_{m_{2n+1}\neq\,m_1,\dots,m_{2n}}\frac{1}{m_{2n+1}}$$
leads to $$\sum_{m_1}\frac{pt_{m_1}^{(0)}}{m_1^{2n-2s+1}}\sum_{m_i}\frac{1}{m_i^{2s}}$$
The full sum will receive from the corrections of its right neighbors $2s-1$ power contributions. We pick $2s$ indices amongst $2n$. The left most sum corresponding to these indices will play the role of full sum. The other $2s-1$ sums will unload onto the latter sum. They can unload directly onto that sum or unload within their own group and so there are $(2s-1)!$ possible ways. The $(2n-2s)$ remaining sums will unload within their own group or on onto the first sum in $m_1$, in $(2n-2s)!$ possible ways.
Then,
$$D_s=(2n-2s)!(2s-1)!\bn{2n}{2s}=\frac{(2n)!}{2s}$$

The integer coefficient $G_k$ counts the number of withdrawal processes of
$$\sum_{m_0}\frac{pt_{m_0}^{(0)}}{m_0}\sum_{m_1\neq m_0}\frac{pt_{m_1}^{(0)}}{m_1}\sum_{m_2\neq m_0,m_1}\frac{1}{m_2}\dots\sum_{m_i\neq m_0,\dots,\,m_{i-1}}\frac{1}{m_i}\dots
\sum_{m_{2n+1}\neq m_0,m_1,\dots, m_{2n}}\frac{1}{m_{2n+1}}$$ leading to
$$\sum_{m_0}\frac{pt_{m_0}^{(0)}}{m_0^{2k+1}}\sum_{m_1}\frac{pt_{m_1}^{(0)}}{m_1^{2n+2-2k-1}}$$
It suffices to pick $2k$ indices amongst the $2n$ right sums which will globally unload onto the first sum in $(2k)!$ possible ways. The remaining $2n-2k$ sums unload onto the second sum in $(2n-2s)!$ possible ways.
Then,
$$G_k=(2n-2k)!(2k)!\bn{2n}{2k}=(2n)!$$
This ends the proof of Lemma $3$.

\subsubsection{Where we resolve the truncated convolutions}

This is the last step of our mod $p^3$ adventure. \\
We first deal with the sum of Corollary $2$ that involves the coefficients $D_k$'s. First, we introduce some new notations.

\newtheorem{Notation}{Notation}
\begin{Notation} Let
\begin{eqnarray}\mcu&:=&\sum_{k=1}^{n-1} \frac{H_{2k}\,S_{p-1-2n+2k}}{2k}\qquad\mpt\\
\mcd&:=&-\sum_{k=1}^{n-1} \frac{H_{2k}\,H_{2n-2k}}{2k}\qquad\qquad\negthickspace\negthickspace\!\!\nts\mpt
\end{eqnarray}
\end{Notation}

\noindent We show that the second sum $\mcd$ is a sum of two truncated convolutions, namely one of divided Bernoulli numbers and one of divided Bernoulli numbers with Bernoulli numbers, as in the following lemma.
\begin{Lemma}
$$\mcd=p^2\Bigg(\sum_{i=p+1-2n}^{p-3}\mb_i\mb_{2(p-1)-2n-i}+\sum_{i=p+1-2n}^{p-3}\frac{B_iB_{2(p-1)-2n-i}}{i}\Bigg)\;\mpt$$
\end{Lemma}
\textsc{Proof of Lemma $4$}. By Result $7$, we know that
\begin{eqnarray*}
H_{2k}&=&\frac{2k}{2k+1}p\,B_{p-1-2k}\qquad\qquad\;\;\,\mpd\\
H_{2n-2k}&=&\frac{2n-2k}{2n-2k+1}p\,B_{p-1-(2n-2k)}\;\mpd
\end{eqnarray*}
Using the latter two congruences in $\mcd$, we get
\begin{equation}
\mcd=-p^2\sum_{k=1}^{n-1}\frac{B_{p-1-2k}}{p-1-2k}.\,\frac{B_{p-1-(2n-2k)}}{p-1-(2n-2k)}.\,(2n-2k)\mpt,
\end{equation}
since the two denominators in $(24)$ can be taken modulo $p$ by von Staudt-Clausen's theorem since $p-1$ does not divide $2k$ nor $2n+2-2k$ when $k$ ranges between integers $1$ and $n-1$. \\
In turn, we write after doing the change of indices $i=p-1-2n+2k$,
\begin{equation}
\mcd=p^2\sum_{i=p+1-2n}^{p-3}(i+1)\frac{B_i}{i}.\,\frac{B_{2(p-1)-2n-i}}{2(p-1)-2n-i}\mpt
\end{equation}
Hence Lemma $4$. The following lemma states a similar result for the sum $\mcu$.
\begin{Lemma}
$$\mcu=-p^2\sum_{i=p+1-2n}^{p-3}\frac{B_iB_{2(p-1)-2n-i}}{i}\mpt$$
\end{Lemma}
\noindent We pleasantly see that some of the terms in Lemma $4$ and $5$ cancel each other. \\\\
\textbf{From now on, whenever this is omitted, it will be understood that we work modulo $p$.} \\\\
\noindent We must compute
$$\sum_{i=p+1-2n}^{p-3}\mb_i\mb_{2(p-1)-2n-i}$$
We will deal with the cases $2n=p-3$ and $2n=p-5$ separately. \\\\Assume for now that $2n<p-5$.\\\\
We will consider the full convolution instead, namely
$$\mathcal{D}:=\sum_{i=2}^{2(p-1)-2n-2}\mb_i\mb_{2(p-1)-2n-i}$$
We decompose it as
\begin{equation}\begin{split}
&\sum_{i=2}^{p-1-2n}\mb_i\mb_{2(p-1)-2n-i}+\sum_{i=p-1}^{2(p-1)-2n-2}\mb_i\mb_{2(p-1)-2n-i}\\
&\\
&+\sum_{i=p+1-2n}^{p-3}\mb_i\mb_{2(p-1)-2n-i}
\end{split}\end{equation}
Moreover, we may group the first two terms together as
\begin{equation}
2\sum_{i=2}^{p-1-2n}\mb_i\mb_{2(p-1)-2n-i}
\end{equation}
Further, for this range of $i$ except for the upper bound, we have $p-1>2n+i$. Hence K\"ummer's congruence applies and yields
$$\mb_{2(p-1)-2n-i}=\mb_{p-1-2n-i}$$
Thus, we get
\begin{equation}\begin{split}
\md\;=\;&2\;\mb_{p-1-2n}\,\mb_{p-1}+2\sum_{i=2}^{p-1-2n-2}\mb_i\,\mb_{p-1-2n-i}\\
&+\sum_{i=p+1-2n}^{p-3}\mb_i\mb_{2(p-1)-2n-i}
\end{split}
\end{equation}
Notice that the first sum of $(28)$ is now a convolution of divided Bernoulli numbers. We will apply Miki's identity to that convolution and to $\md$ as well. In the latter case (take $n:=2(p-1)-2n$ in Result $10$), we get:
\begin{equation}\begin{split}
\md\;=\;&\sum_{i=2}^{2(p-1)-2n-2}\binom{2(p-1)-2n}{i}\mb_i\mb_{2(p-1)-2n-i}\\&\\&+\,2\,\mh_{2(p-1)-2n}\,\mb_{2(p-1)-2n}
\end{split}\end{equation}
When $p-1-2n\leq i\leq p-1$, we have
$2(p-1)-2n-i+1\leq p$, so that
$$\binom{2(p-1)-2n}{i}=\frac{(2(p-1)-2n)\dots (2(p-1)-2n-i+1)}{i!}=0$$
However, when $i=p-1$ or $i=p-1-2n$, the product of divided Bernoulli numbers reads $\mb_{p-1-2n}\mb_{p-1}$ and $p$ divides the denominator of $B_{p-1}$ by Von Staudt-Clausen's theorem. In that case, the binomial coefficient should not get canceled modulo $p\mathbb{Z}_p$. \\
Also, when $i\leq p-1-2n-2$, we have $p-1-2n-i>0$. Hence, we get after applying K\"ummer's congruence,
\begin{equation}\begin{split}
\md\;=\;&2\,\sum_{i=2}^{p-1-2n-2}\binom{2(p-1)-2n}{i}\mb_i\mb_{p-1-2n-i}\\&\\
&+2\binom{2(p-1)-2n}{p-1}\mb_{p-1-2n}\mb_{p-1}\\&\\&+\,2\,\mh_{2(p-1)-2n}\,\mb_{2(p-1)-2n}
\end{split}\end{equation}
Another application of Miki's identity yields:
\begin{equation}\begin{split}
2\sum_{i=2}^{p-1-2n-2}\mb_i\,\mb_{p-1-2n-i}\,=\,&2\,\sum_{i=2}^{(p-1-2n)-2}\binom{p-1-2n}{i}\mb_i\,\mb_{p-1-2n-i}\\&\\&+4\,\mh_{p-1-2n}\,\mb_{p-1-2n}
\end{split}\end{equation}
By gathering all the information so far, we obtain:
\begin{equation}
\begin{split}
&\sum_{i=p+1-2n}^{p-3}\mb_i\mb_{2(p-1)-2n-i}\,=\\&2\,\sum_{i=2}^{p-1-2n-2}\Big\lbrace\binom{2(p-1)-2n}{i}-\binom{p-1-2n}{i}\Big\rbrace\mb_i\mb_{p-1-2n-i}\\
&+2\Bigg\lbrace\binom{2(p-1)-2n}{p-1}-1\Bigg\rbrace\mb_{p-1-2n}\mb_{p-1}\\
&+2\mb_{2(p-1)-2n}\,\mh_{2(p-1)-2n}-4\mb_{p-1-2n}\,\mh_{p-1-2n}
\end{split}
\end{equation}
Next, we apply the Chu-Vandermonde convolution (Result $12$ with $m:=p-1$, $n:=p-1-2n$ and $r:=i$),
\begin{equation}
\binom{2(p-1)-2n}{i}=\sum_{j=0}^{i}\bn{p-1}{j}\bn{p-1-2n}{i-j}
\end{equation}
Applying further Lemma $5$ of \cite{LEV2}, we derive
\begin{equation}
\binom{2(p-1)-2n}{i}=\sum_{j=0}^{i}(-1)^j\bn{p-1-2n}{i-j}
\end{equation}
Denoting the first sum to the right hand side of $(32)$ by $\mathcal{S}$, we get successively,
\begin{eqnarray}
\mathcal{S}&=&2\sum_{i=2}^{p-1-2n-2}\Bigg\lbrace\sum_{j=1}^i(-1)^j\bn{p-1-2n}{i-j}\Bigg\rbrace\mb_i\mb_{p-1-2n-i}\\
&=&2\sum_{i=2}^{p-1-2n-2}\Bigg\lbrace\sum_{j=0}^{i-1}(-1)^j\bn{p-1-2n}{j}\Bigg\rbrace\mb_i\mb_{p-1-2n-i}\\
&=&2\sum_{i=2}^{p-1-2n-2}(-1)^{i-1}\bn{p-2-2n}{i-1}\mb_i\mb_{p-1-2n-i}\\
&=&-2\sum_{i=2}^{p-1-2n-2}\bn{p-2-2n}{i-1}\mb_i\mb_{p-1-2n-i}
\end{eqnarray}
For the calculation of the alternating sum of binomial coefficients in $(36)$, see for instance \cite{BO}.\\
We now consider the sum of $(38)$ without its factor. We distinguish between two cases. Either we sum an even number of terms or we sum an odd number of terms. In the latter case, the middle index corresponds to
$$i=\frac{p-1}{2}-n$$
and we have $$\bn{p-2-2n}{\frac{p-1}{2}-n}=\bn{p-2-2n}{\frac{p-1}{2}-n-1}$$
Then, if we consider twice the sum of $(38)$, the middle term contributes to
\begin{equation}
\bn{p-1-2n}{\frac{p-1}{2}-n}\;\mb_{\frac{p-1}{2}-n}\mb_{\frac{p-1}{2}-n}
\end{equation}
And the other terms contribute to
\begin{equation}\sum_{\begin{array}{l}2\leq i\leq p-1-2n-2\\i\neq\,\frac{p-1}{2}-n\end{array}}\bn{p-1-2n}{i}\mb_i\mb_{p-1-2n-i}
\end{equation}
If we sum an even number of terms, we get directly
\begin{equation}\sum_{i=2}^{p-1-2n-2}\bn{p-1-2n}{i}\mb_i\mb_{p-1-2n-i}\end{equation}
Therefore, after a new application of Miki's identity, we get:
\begin{equation}
\mathcal{S}\;=\;2\;\mh_{p-1-2n}\;\mb_{p-1-2n}\;-\sum_{i=2}^{p-1-2n-2}\mb_i\,\mb_{p-1-2n-i}
\end{equation}
It remains to deal with the second row of Congruence $(32)$. Denote this term by $\mt$. First and foremost, we also apply the Vandermonde-Chu equality and obtain
\begin{equation}\binom{2(p-1)-2n}{p-1-2n}=\sum_{j=0}^{p-1-2n}\binom{p-1}{j}\binom{p-1-2n}{p-1-2n-j}\end{equation}
Second, by Proposition $6$ of \cite{LEV2}, we have (with $\mh_0=0$):
\begin{equation}\binom{p-1}{j}=(-1)^j(1-p\mh_j)\mpd\end{equation}
We thus have,
\begin{equation}\frac{1}{2}p^2\mt=-(p+1)p\,B_{p-1}p\,\mb_{p-1-2n}\sum_{j=1}^{p-1-2n}(-1)^j(1-p\mh_j)\binom{p-1-2n}{p-1-2n-j}\mpt\end{equation}
Denote the binomial sum above by $BS$. We have
\begin{eqnarray}BS&=&\sum_{j=0}^{p-2-2n}(-1)^j\binom{p-1-2n}{j}-p\sum_{j=1}^{p-1-2n}(-1)^j\binom{p-1-2n}{j}\mh_j\\
&=&-1+\frac{p}{p-1-2n}\end{eqnarray}
The alternating binomial harmonic sum identity is due to Spivey and appears in \cite{SP}, identity $20$.
The sum $BS$ must be considered modulo $p^2$, hence we get:
\begin{equation}
\frac{1}{2}p^2\mt=-(p+1)p\,B_{p-1}p\,\mb_{p-1-2n}\Big(-1-\frac{p}{2n+1}\Big)
\end{equation}

We are now ready to give the almost final form of Congruence $(32)$. We state it in the following proposition. We used the fact that $\mh_{p-1-2n}$ is a $p$-adic integer, hence
$$\mb_{p-1-2n}\mh_{p-1-2n}=\mb_{2(p-1)-2n}\mh_{p-1-2n}\qquad\text{mod $p\,\mathbb{Z}_p$}$$
We also list below the special cases $2n\in\lbrace p-3,p-5\rbrace$ which can be easily processed by hand.

\newtheorem{Proposition}{Proposition}
\begin{Proposition} \textbf{Assume $4\leq 2n\leq p-7$.}
\begin{equation}\begin{split}
p^2\sum_{i=p+1-2n}^{p-3}\mb_i\mb_{2(p-1)-2n-i}
=&-p^2\sum_{i=2}^{p-1-2n-2}\mb_i\mb_{p-1-2n-i}\\&
+2p\,\mb_{2(p-1)-2n}\,p\big(\mh_{2(p-1)-2n}-\mh_{p-1-2n}\big)\\
&+2(p+1)\,p\,B_{p-1}\Big(1+\frac{p}{2n+1}\Big)\,p\,\mb_{p-1-2n}\mpt
\end{split}
\end{equation}

\noindent\textbf{Case $2n=p-3$.}
\begin{equation}
p^2\sum_{i=4}^{p-3}\mb_i\mb_{p+1-i}=2p\,\mb_{p+1}+p^2\mb_2+2p\mb_2(p\,B_{p-1})\mpt
\end{equation}
\textbf{Case $2n=p-5$.}
\begin{equation}
p^2\sum_{i=6}^{p-3}\mb_i\mb_{p+3-i}=\frac{7}{720}\,p^2+2p\,\mb_{p+3}+2p\,\mb_4(p\,B_{p-1})\mpt
\end{equation}
\end{Proposition}

\begin{Remark} For a list of the first Bernoulli numbers, see \cite{BE} or \cite{NI}.\\
Congruence $(51)$ got tested on prime $p=11$. Both sides of the congruence evaluated successfully to $1210$ modulo $11^3$. \\
Congruence $(50)$ got verified on prime $p=17$ leading to $1734$ modulo $17^3$ on both sides.
\end{Remark}

\noindent Before moving further, it is time to step back and go back to the expression modulo $p^3$ for the Stirling numbers as provided in \cite{LEV}. This is Result $8$ of the introduction.
First, we introduce some notations for convolutions of divided Bernoulli numbers.

\begin{Notation} Let
$$\mc\mb(n):=\sum_{i=2}^{n-2}\mb_i\mb_{n-i}$$
\end{Notation}

By Result $8$ applied with $2k=p-1-2n$, we have when $0\leq 2n\leq p-5$,
\begin{equation}
\left[\begin{array}{l}\;p\\2n+1\end{array}\right]=-\;\frac{1}{p-1-2n}\Bigg(p\,B_{p-1-2n}-p^2\,\sum_{i=2}^{p-3-2n}\mb_i\,B_{p-1-2n-i}\Bigg)\mpt
\end{equation}
Moreover, notice that twice the sum participating into $(52)$ is congruent modulo $p\,\mathbb{Z}_p$ to
$$-(2n+1)\sum_{i=2}^{p-3-2n}\mb_i\mb_{p-1-2n-i}$$
Hence we get the following result:
\begin{Proposition} Assume $0\leq 2n\leq p-5$. Then, we have:
\begin{equation}
\frac{p^2}{2}\mc\mb(p-1-2n)=\left[\begin{array}{l}\;p\\2n+1\end{array}\right]\;+\,p\,\mb_{p-1-2n}\qquad\mpt
\end{equation}
In particular, we have
\begin{equation}
\frac{p^2}{2}\mc\mb(p-1)=p\,\mb_{2(p-1)}-\frac{1}{2}p^2\mb_{p-1}^2\qquad\qquad\qquad\mpt
\end{equation}
\end{Proposition}
\noindent Congruence $(54)$ is Congruence $(53)$ using the original fact from \cite{SU2} that
\begin{equation}(p-1)!=-p\,\mb_{p-1}+p\mb_{2p-2}-\frac{1}{2}p^2\mb_{p-1}^2\qquad\mpt\end{equation}
See also Corollary $6$ of \cite{LEV} where Sun's formula gets proven using part of the methodology of the current paper. \\

We now process the last sum of Corollary $2$. Again, we split the sum, this time into four contributing parts. We introduce some notations in order to refer to the different parts.
\begin{Notation}
\begin{eqnarray}
\mct&:=&\sum_{k=1}^{n-1} S_{p-1-2k}S_{p-1-2n+2k}\;\;\;\;\mpt\\
\mcq&:=&-\sum_{k=1}^{n-1} S_{p-1-2k}H_{2n-2k}\qquad\mpt\\
\mcf&:=&-\sum_{k=1}^{n-1} H_{2k}S_{p-1-2n+2k}\qquad\mpt\\
\mcs&:=&\sum_{k=1}^{n-1} H_{2k}H_{2n-2k}\qquad\qquad\;\;\mpt\\
\end{eqnarray}
\end{Notation}

We see with a simple change of indices that $\mcq$ and $\mcf$ are equal. Moreover, their common value modulo $p^3$ is
$$-p^2\sum_{k=1}^{n-1} \frac{2k}{2k+1}\,B_{p-1-2k}B_{p-1-2n+2k}$$
We can group the sums again adequately and obtain:
\begin{equation}\begin{split}
\sum_{i=3}^6\mc_i&=p^2\sum_{k=1}^{n-1} \Bigg(1-\frac{2k}{2k+1}\Bigg)B_{p-1-2k}B_{p-1-2n+2k}\\
&-p^2\sum_{k=1}^{n-1} \frac{2k}{2k+1}\Bigg(1-\frac{2n-2k}{2n-2k+1}\Bigg)B_{p-1-2k}B_{p-1-2n+2k}\\&\qquad\qquad\qquad
\qquad\qquad\qquad\qquad\qquad\qquad\qquad\qquad\qquad\mpt
\end{split}\end{equation}
Hence, we have:
\begin{equation}\begin{split}
\sum_{i=3}^6\mc_i&=-p^2\sum_{k=1}^{n-1} \mb_{p-1-2k}B_{p-1-2n+2k}\\&+p^2\sum_{k=1}^{n-1} \frac{2k}{2k+1}B_{p-1-2k}\mb_{p-1-2n+2k}\qquad\mpt
\end{split}\end{equation}
So finally,
\begin{equation}\begin{split}
\sum_{i=3}^6\mc_i&=-p^2\sum_{k=1}^{n-1} \mb_{p-1-2k}B_{p-1-2n+2k}\\&+p^2\sum_{k=1}^{n-1} B_{p-1-2k}\mb_{p-1-2n+2k}\\&
+p^2\sum_{k=1}^{n-1} \mb_{p-1-2k}\mb_{p-1-2n+2k}\qquad\mpt
\end{split}\end{equation}
And after cancelation of the first two terms, we get
\begin{eqnarray}
\sum_{i=3}^6\mc_i&=&p^2\sum_{k=1}^{n-1} \mb_{p-1-2k}\mb_{p-1-2n+2k}\\
&=&p^2\sum_{i=p+1-2n}^{p-3}\mb_i\mb_{2(p-1)-2n-i}\qquad\mpt
\end{eqnarray}
\subsubsection{Final expression and confrontation with the result issued from Manner I}
By Sun's result taken from \cite{SU2} and expressed in Result $4$, $(a)$ of $\S\,1.2$, we have
\begin{equation}
\text{When $2n\neq p-3$,}\;\;\;H_{2n}=2n\Big(p\,\mb_{2p-2-2n}-2p\,\mb_{p-1-2n}\Big)\mpt
\end{equation}
Also, we have:
\begin{eqnarray}
S_{p-1-2n}&=&(p-1-2n)\,p\,\mb_{p-1-2n}\;\;\;\;\,\,\mpt\\
S_{2p-2n-2}&=&(2p-2n-2)\,p\,\mb_{2p-2n-2}\;\;\mpt\\
p^2\mb_{2p-2n-2}&=&p^2\mb_{p-1-2n}\qquad\qquad\qquad\;\;\,\mpt
\end{eqnarray}
Using these facts as well as Corollary $2$ and the expressions for the sums $\mc_i$'s with $1\leq i\leq 6$, we have
\begin{equation}\begin{split}
\text{When $2n\neq p-3$,}\;\;\;\left[\begin{array}{l}\;p\;\\2n+1\end{array}\right]=&-\frac{p^2}{2}\sum_{i=p+1-2n}^{p-3}\mb_i\mb_{2(p-1)-2n-i}\\
&+p\,\mb_{2p-2n-2}\\
&+\Big(p(p\,B_{p-1})_1-p-2\Big)\,p\,\mb_{p-1-2n}\\&\qquad\qquad\qquad\qquad\qquad\qquad\mpt
\end{split}\end{equation}
When $2n=p-3$, we use the Sun congruences (cf Results $2$ and $4$ $(b)$),
\begin{eqnarray}
H_{p-3}&=&\big(\frac{1}{2}-3\,B_{p+1}\big)p-\frac{4}{3}p^2\;\;\;\mpt\\
S_2&=&\,p\,B_2+p^2\,B_1\qquad\qquad\!\!\;\;\;\mpt\\
S_{p+1}&=&\,p\,B_{p+1}\qquad\qquad\qquad\;\,\;\;\mpt
\end{eqnarray}
And so we have,
\begin{equation}
\left[\begin{array}{l}\;\;p\\p-2\end{array}\right]=-\frac{p^2}{2}\sum_{i=4}^{p-3}\mb_i\mb_{p+1-i}+p\,B_{p+1}+\frac{p^2}{12}(p\,B_{p-1})_1+\frac{p^2}{3}-\frac{p}{6}\;\;\mpt
\end{equation}
Before we use Proposition $1$, we need a lemma.
\begin{Lemma} Assume $2n<p-5$. The following congruence holds.
\begin{equation}
p\bigg(\mh_{2(p-1)-2n}-\mh_{p-1-2n}\bigg)=1+\frac{p}{2n+1}\mpd
\end{equation}
\end{Lemma}
\textsc{Proof of Lemma $6$}. We have
\begin{eqnarray}
p\bigg(\mh_{2(p-1)-2n}-\mh_{p-1-2n}\bigg)&=&p\Bigg(\frac{1}{p-1-2n+1}+\frac{1}{p-1-2n+2}+\dots+\frac{1}{p-1-2n+2n}\notag\\&&+\frac{1}{p-1-2n+2n+1}\notag\\
&&+\frac{1}{p-1-2n+2n+2}+\dots+\frac{1}{p-1-2n+2n+p-1-2n}\Bigg)\notag\\&&\notag\\
&=&-p\,\mh_{2n}+1+p\,\mh_{p-(2n+2)}\qquad\mpd\\
&=&-p\,\mh_{2n}+1+p\,\mh_{2n+1}\qquad\qquad\!\!\mpd\\
&=&1+\frac{p}{2n+1}\qquad\qquad\qquad\qquad\;\;\mpd
\end{eqnarray}
with $(77)$ following from an application of Wolstenholme's theorem (see \cite{LEV2}, $\S\,4$ p. $22$).
Assume $2n<p-5$. In light of Lemma $6$, Proposition $1$ rewrites as
\begin{equation}\begin{split}
p^2\sum_{i=p+1-2n}^{p-3}\mb_i\mb_{2(p-1)-2n-i}=&-p^2\mc\mb(p-1-2n)\\&+2p\,\mb_{2(p-1)-2n}-2(p+p^2)\,\mb_{p-1-2n}\\&+2(p\,B_{p-1})_1\,p^2\mb_{p-1-2n}\;\;\mpt
\end{split}\end{equation}
\noindent Point $(i)$ of Theorem $1$ is a rearrangement of Congruence $(79)$ and is thus proven. Points $(ii)$ and $(iii)$ are respectively Congruences $(50)$ and $(51)$ of Proposition $1$.\\
The conjunction of Congruences $(70)$ and $(79)$ leads precisely to Congruence $(53)$ of Proposition $2$, which is the result arising from Manner I.\\
It is also a straightforward (but fastidious) verification that the formulas from Manner $1$ and Manner $2$ for the Stirling number $\left[\begin{array}{l}
\;\;p\\p-4\end{array}\right]$ modulo $p^3$ coincide. And the formula arising from Manner $2$ for computing $\left[\begin{array}{l}
\;\;p\\p-2\end{array}\right]$ yields the same result as the one issued from a direct calculation (see last row of Result $8$ of the introduction). \\

The next part deals with multiple harmonic sums and relates them in particular modulo $p^3$ to the Stirling numbers in a reverse way than the trivial usual way, that is $A_{p-1-j}=(p-1)!\,\as_j$.

\section{The multiple harmonic sums modulo $p^3$}

This part is based on the use of Newton's formulas which we recall right below.
\begin{Result} (Newton's formula, see \cite{JAC}). Let $x_1,\dots,x_m$ be complex numbers.
Let $\mcp_k=x_1^k+\dots x_m^k$ and $\ma_k=\sum_{1\leq i_1<\dots<i_k\leq m}x_{i_1}\dots x_{i_k}$. \\
Then, for $k=0,1,\dots,m$ we have:
\begin{equation}\mcp_k-\ma_1\mcp_{k-1}+\ma_2\mcp_{k-2}+\dots+(-1)^{k-1}\ma_{k-1}\mcp_1+(-1)^{k}k\ma_k=0\end{equation}
\end{Result}
\noindent Our complex numbers are chosen to be the first $p-1$ reciprocals of integers. And so $\mcp_k=H_k$ and $\ma_k=\as_k$.
Then, Equality $(80)$ rewrites as:
\begin{equation}
\as_k=\frac{(-1)^{k-1}}{k}\Bigg(H_k+\sum_{r=1}^{k-1}(-1)^r\as_r\,H_{k-r}\Bigg)\qquad\forall k=1,\dots,p-1
\end{equation}
Note, when $k=1$, $\as_1=H_1$, hence the equality still holds. \\
We will follow Sun's idea exploited in a different context. By Bayat's result, we know that $p$ divides $H_1,\dots,H_{p-3}$. Moreover, by Sun's result $4$ $(c)$ or by Corollary $4$ of \cite{LEV}, we know that $p$ also divides $H_{p-2}$. Then we see from $(81)$ that $p$ divides all of $\as_1,\dots,\as_{p-1}$ (note the latter fact can also be seen directly from $A_{p-1-j}=(p-1)!\,A_j^{\star}$ and the fact that $p$ divides the Stirling numbers $A_1,\dots,A_{p-2}$. This implies in turn that $p^2$ divides the sum of $(81)$ and so we have,
\begin{equation}
\as_k=\frac{(-1)^{k-1}}{k}\,H_k\mpd\qquad\qquad\qquad\;\;\;\; \forall\,k=1,\dots,p-1
\end{equation}
Now if we work modulo $p^3$, the multiple harmonic sums must be considered modulo $p^2$ in the sum of $(81)$. Thus, we get
\begin{equation}
\as_k=\frac{(-1)^{k-1}}{k}\Bigg(H_k-\sum_{r=1}^{k-1}\frac{H_r\,H_{k-r}}{r}\Bigg)\mpt\qquad\forall k=1,\dots,p-1
\end{equation}
Moreover, by Corollary $5.1$ of \cite{SU2}, we have
$$H_r=\frac{r}{r+1}p\,B_{p-1-r}\mpd\qquad\text{when $r<p-1$}$$
Then,
\begin{equation}
\as_k=\frac{(-1)^{k-1}}{k}\Bigg(H_k-p^2\,\sum_{r=1}^{k-1}(k-r)\mb_{p-1-r}\mb_{p-1-k+r}\Bigg)\mpt\qquad\forall k=1,\dots,p-1
\end{equation}
If in the sum of $(84)$, we do the change of indices $i=p-1-k+r$, we obtain:
\begin{equation}
\as_k=\frac{(-1)^{k-1}}{k}\Bigg(H_k+p^2\,\sum_{i=p-k}^{p-2}(i+1)\mb_i\mb_{2(p-1)-k-i}\Bigg)\mpt\qquad\forall k=1,\dots,p-1
\end{equation}
It follows that
\begin{equation}
\as_k=\frac{(-1)^{k-1}}{k}\Bigg(H_k-\frac{k}{2}p^2\sum_{i=p+1-k}^{p-3}\mb_i\mb_{2(p-1)-k-i}\Bigg)\mpt\qquad\forall k=1,\dots,p-1
\end{equation}
When $k$ is odd, $2(p-1)-k-i$ and $i$ have distinct parity, hence the sum of $(86)$ vanishes. If $k$ is even, say $k=2n$, the sum is the well studied truncated sum of $\S\,2.2.3$, namely $\mt\mc\mb(p-1-2n,p-3)$ using the notations from the introduction.\\

\noindent From there, point (I) of Theorem $2$ can be derived from our Theorem $1$ and from Sun's result listed as Result $4$ of the introduction.\\
As for point (II) of Theorem $2$, it can be obtained directly from $(82)$ by multiplying each side of the congruence by $p\,w_p$. \\
Finally, point (III) follows from the fact that $p$ divides $H_1,\dots,H_{p-2}$ and so $p$ divides $\as_1,\dots,\as_{p-2}$ by $(82)$, as we have already seen before. Hence,
it follows from Glaisher's result for $(p-1)!$ modulo $p^2$ that
\begin{equation}A_{p-1-j}=(p-1)!\as_j=(p\,B_{p-1}-p)\,\as_j\mpt\qquad\forall j=1,\dots,p-2\end{equation}

\noindent It is easy to check that $A_1=(p\,B_{p-1}-p)\as_{p-2}$ and $A_2=(p\,B_{p-1}-p)\as_{p-3}$, this modulo $p^3$.
We note moreover that the conjunction of Proposition $2$ of $\S\,2.2.3$ and of congruences $(86)$ and $(87)$ provides a third independent proof for the relationship between the truncated convolution and the full convolution of interests modulo $p^3$.\\

We end this part by proving Theorems $0.0.,0.1.,0.2.$ and $0.3.$ \\\\By Zhao's congruence $(14)$ in \cite{ZO}, we have
\begin{equation}
H_1=p\as_2\;\;\text{mod}\,p^4
\end{equation}
By Equality $(81)$ with $k=2$, we have
\begin{eqnarray}
\as_2&=&-\frac{1}{2}(H_2-\as_1H_1)\\
&=&-\frac{1}{2}(H_2-H_1^2)
\end{eqnarray}
By Wolstenholme's theorem, $p^2$ divides $H_1$, hence $p^4$ divides $H_1^2$. Thus, we get
\begin{eqnarray}
\as_2&=&-\frac{1}{2}H_2\qquad\qquad\;\;\;\;\;\;\mpt\\
&=&-p(\mb_{2p-4}-2\mb_{p-3})\mpt
\end{eqnarray}
Then, $$H_1=-p^2(\mb_{2p-4}-2\mb_{p-3})\;\;\,\mpq$$
This settles Theorem $0.0$. \\
Regarding Theorem $0.1$, it can be derived from Theorem $0.0$ in the following way. 
\begin{eqnarray}
A_{p-3}&=&(p\,B_{p-1}-p)\as_2\qquad\qquad\qquad\qquad\mpt\\
&=&(p\,B_{p-1}-p)\frac{H_1}{p}\qquad\,\;\;\;\;\;\qquad\qquad\mpt\\
&=&(p\,B_{p-1}-p)\,p\,\mb_{p-3}-\frac{p^2}{2}\sum_{a=1}^{p-1}\frac{q_a^2}{a^2}\,\;\;\mpt
\end{eqnarray}
As far as Theorem $0.2$, we derive from Proposition $2$ applied with $2n=2$ and from Theorem $0.1$ that
$$\frac{p^2}{2}\,\mc\mb(p-3)=(p\,B_{p-1}-p+1)p\,\mb_{p-3}-\frac{p^2}{2}\sum_{a=1}^{p-1}\frac{q_a^2}{a^2}\qquad\mpt$$
Hence,
$$p\,\mc\mb(p-3)=2(p\,B_{p-1}-p+1)\,\mb_{p-3}-p\,\sum_{a=1}^{p-1}\frac{q_a^2}{a^2}\qquad\,\,\,\mpd$$
Thus, in terms of the Agoh-Giuga quotient and then the Wilson quotient, we obtain:
\begin{eqnarray}
\mc\mb(p-3)&=&2\Bigg(\Big(\frac{p\,B_{p-1}+1}{p}\Big)-1\Bigg)\mb_{p-3}-\sum_{a=1}^{p-1}\frac{q_a^2}{a^2}\qquad\mpu\\
&=&2\,w_p\,\mb_{p-3}-\sum_{a=1}^{p-1}\frac{q_a^2}{a^2}\qquad\qquad\qquad\qquad\qquad\mpu
\end{eqnarray}
Regarding the convolution of order $p-5$, Theorem $0.3.$ is derived from a special case of Theorem $1$ point $(i)$ applied with $2n=4$. It provides $\mc\mb(p-5)$ modulo $p$. We then obtain the Stirling numbers on $5$ disjoint cycles modulo $p^3$ by the congruence $(53)$ of Proposition $2$. The multiple harmonic sum $\as_4$ modulo $p^3$ is finally deduced from the congruence of Theorem $2$, point (II), or directly from Congruence $(86)$. \\
We have thus "resolved" the convolutions $\mc\mb(p-1), \mc\mb(p-3)$ and $\mc\mb(p-5)$ to the modulus $p$.
\begin{Remark}
In \cite{ZO2} p. $97$, Jianqiang Zhao proved using generalized multiple harmonic sums that
\begin{eqnarray}
\forall p\geq 7,\;\;CB(p-3)&=&-2\,B_{p-3}\qquad\qquad\;\;\;\mpu\\
\forall p\geq 9,\;\;CB(p-5)&=&-\frac{2}{3}\,B_{p-3}^2-2\,B_{p-5}\;\;\mpu
\end{eqnarray}
where $CB$ denotes convolutions of ordinary Bernoulli numbers.
\end{Remark}

\section{Concluding words}
The $p$-adic analysis on the polynomial $X^{p-1}+(p-1)!\in\mathbb{Z}_p[X]$ joint with some century(ies) old mathematics allows to bypass more recent and technical results concerning Miki's identity or the knowledge of the generalized harmonic numbers modulo $p^3$ which uses in particular a generalization of K\"ummer's congruences to the modulus $p^2$ by Zhi-Hong Sun. \\
We have resolved the convolutions only for special cases and there is far more work to be done. \\
At the level of the modulus $p^4$, it is to expect that the whole set of tools detailed here will prove useful and non intersecting. \\\\
\indent \textit{Email address:} clairelevaillant@yahoo.fr

\end{document}